%% file: ap-posetout.tex
\def\hldest#1#2#3{}
\input fontmac
\input mathmac
\input epsf
\input eplain

\def\divides{\backslash}

\def\rank{\op{\rm rank}}
\def\lcm{\op{\rm lcm}}
\def\ref#1{\special{ps:[/pdfm { /big_fat_array exch def big_fat_array 1 get 0
0 put big_fat_array 1 get 1 0 put big_fat_array 1 get 2 0 put big_fat_array pdfmnew } def}%
[\hlstart{name}{}{bib#1}#1\hlend]}
\def\hat{\widehat}
\def\tilde{\widetilde}
\def\bar{\overline}
\def\filter{\op{$\uparrow$}}
\def\ideal{\op{$\downarrow$}}

\magnification=\magstephalf
\hoffset=40pt \voffset=28pt
\hsize=29pc  \vsize=45pc  \maxdepth=2.2pt  \parindent=19pt
\nopagenumbers
\def\leftheadline{{\rm\folio}\hfil{\eightpoint GOH, HAMDAN, AND SAKS}\hfil}
\def\rightheadline{\hfil{\eightpoint THE LATTICE OF ARITHMETIC PROGRESSIONS}\hfil{\rm\folio}}
\headline={\ifodd\pageno{\ifnum\pageno<2\hfil\else\rightheadline\fi}\else\leftheadline\fi}

\enablehyperlinks

\ifpdf
  \hlopts{bwidth=0}
  \pdfoutline goto name {intro} {Introduction}%
  \pdfoutline goto name {numprogs} {The number of arithmetic progressions}%
  \pdfoutline goto name {chains} {Chains and the order complex}%
  \pdfoutline goto name {coatoms} {Coatoms}%
  \pdfoutline goto name {homology} {Homology groups of the order complex}%
  \pdfoutline goto name {comodernism} {Left-modularity and comodernism}%
  \pdfoutline goto name {homotopy} {EL-labelability, homotopy type, and complements}%
  \pdfoutline goto name {acks} {Acknowledgements}%
  \pdfoutline goto name {refs} {References}%
\else
  \special{ps:[/PageMode /UseOutlines /DOCVIEW pdfmark}%
  \special{ps:[/Count -0 /Dest (intro) cvn /Title (Introduction) /OUT pdfmark}%
  \special{ps:[/Count -0 /Dest (numprogs) cvn /Title (The number of arithmetic progressions) /OUT pdfmark}%
  \special{ps:[/Count -0 /Dest (chains) cvn /Title (Chains and the order complex) /OUT pdfmark}%
  \special{ps:[/Count -0 /Dest (coatoms) cvn /Title (Coatoms) /OUT pdfmark}%
  \special{ps:[/Count -0 /Dest (homology) cvn /Title (Homology groups of the order complex) /OUT pdfmark}%
  \special{ps:[/Count -0 /Dest (comodernism) cvn /Title (Left-modularity and comodernism) /OUT pdfmark}%
  \special{ps:[/Count -0 /Dest (homotopy) cvn /Title
  (EL-labelability, homotopy type, and complements) /OUT pdfmark}%
  \special{ps:[/Count -0 /Dest (acks) cvn /Title (Acknowledgements) /OUT pdfmark}%
  \special{ps:[/Count -0 /Dest (refs) cvn /Title (References) /OUT pdfmark}%
\fi

\maketitle{The
lattice of arithmetic progressions}{}{Marcel K. Goh, Jad Hamdan, {\rm and} Jonah Saks}
{{\sl Department of Mathematics and Statistics, McGill University}}

\floattext5 \ninebf Abstract.
\ninepoint
This paper concerns the lattice $L_n$ of subsets of $\{1,\ldots,n\}$ that are arithmetic
progressions, under the
inclusion order. For $n\geq 4$, this poset is not graded and thus not semimodular.
We give three independent proofs of the fact that for $n\geq 2$,
$\mu_n(L_n) = \mu(n-1)$, where $\mu_n$ is the M\"obius function of $L_n$ and
$\mu$ is the classical (number-theoretic) M\"obius function. We also show that $L_n$ is comodernistic,
which implies that $L_n$ is EL-labelable. Comodernism is then used to prove
that the order complex $\Delta_n$ of the lattice is either contractible or homotopy equivalent to a sphere.
\smallskip
\noindent\boldlabel Keywords. Lattices, arithmetic progressions, M\"obius function, order complex.

\advsect Introduction
\hldest{xyz}{}{intro}

{\tensc The additive structure} of certain subsets has long been a topic of interest in number
theory and combinatorics.
A class of sets with a great deal of additive structure is the set of {\it arithmetic
progressions}. These are sets of the form
$$\big\{a, a+r, \ldots, a+(k-1)r\big\}$$
where the {\it base point} $a$ and {\it step size} (or simply {\it step})
$r$ are elements of an additive group and the {\it length}
$k$ is an integer. In this paper we take our underlying additive group to be the integers $\ZZ$.
The business of finding arithmetic progressions in sets of integers goes back to a classical 1927 theorem
of B.~L.~van der Waerden~\ref{23},
which states that any colouring of the integers with finitely many colours gives rise
to monochromatic arithmetic progressions of arbitrary length. This was generalised by E.~Szemer\'edi, who
in 1975 proved the existence of arithmetic progressions of arbitrary length in any set of positive upper
density~\ref{22}. More recently, B.~Green and T.~Tao showed that the same conclusion
holds in the primes~\ref{14}.

Set systems consisting of arithmetic progressions have received some attention in the realm of topology.
The topology on $\ZZ$ generated by (infinite) arithmetic progressions $a+k\ZZ$ was used by H.~Furstenberg
to give an alternative proof of the infinitude of primes~\ref{10}. This topology came to be
known as Golomb's topology, after S.~Golomb who studied its properties more systematically in a 1959
paper~\ref{12}.
We will restrict ourselves to a finite subset of $\ZZ$ and study the set of arithmetic progressions itself,
rather than the topology
it forms a basis of. As with any set of subsets, it is partially ordered
by inclusion, and in the present paper we investigate the structure induced by this ordering.

We shall also investigate topological properties of the order complex associated to this lattice.
Several other simplicial
complexes related to number-theoretic objects have recently appeared in the literature.
The simplicial complex of squarefree positive
integers less than or equal to $n$ was studied in a 2011 paper by A.~Bj\"orner~\ref{4},
and a 2017 paper~\ref{7} of
R.~Ehrenborg, L.~Govindaiah, P.~S.~Park, and M.~Readdy introduces a simplicial complex called the van der Waerden
complex ${\rm vdW}(n,k)$, whose facets correspond to arithmetic progressions of length $k$ in
$\{1,\ldots,n\}$. A subsequent paper of B.~Hooper and A.~Van~Tuyl characterised the pairs $(n,k)$ for which
${\rm vdW}(n,k)$ is shellable~\ref{16}.
The simplicial complexes arising from our posets are different in that the vertices are themselves arithmetic
progressions.

\newcount\smallcases
\smallcases=\figcount
Let $[n]$ denote the set $\{1,2,\ldots,n\}$. For $n\geq 1$ we
let $L_n$ denote the partially-ordered set (poset)
of all finite integer arithmetic progressions contained in $[n]$ including trivial progressions of length
$1$ and $2$ as well as the empty set $\emptyset$. When it is convenient, we artificially
define $L_0 = \{\emptyset\}$. Small examples are depicted in Fig.~{\the\smallcases}.
\midinsert
$$\epsfbox{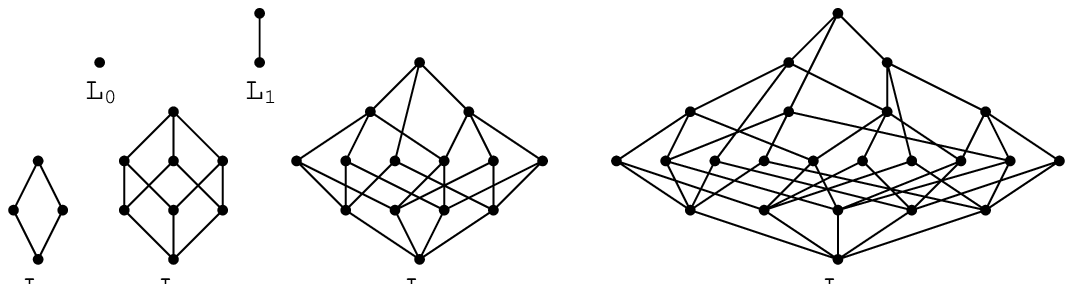}$$
\vskip5pt
\caption{Hasse diagrams of $L_n$ for small values of $n$.}
\endinsert
\goodbreak
The notation $L_n$ is motivated by the fact that $L_n$ is a lattice. The meet of two elements is simply
the set-theoretic intersection, since the intersection of two integer
arithmetic progressions is a (possibly empty)
arithmetic progression. By induction, one finds that the meet of any finite number of points
is well-defined in $L_n$, and this as well as the existence of
a maximum element $12\cdots n$ implies the existence of a join
of two arbitrary elements $x_1, x_2\in L_n$. (Because
$S = \{x \in L_n : x_1\cup x_2\subseteq x\}$ is nonempty (it contains at least $12\cdots n$),
we may set $x_1 \vee x_2 = \bigwedge_{x\in S} x$.)

The poset $L_n$ is not graded for $n\geq 4$. To see this, note that
$$1<14<1234<12345<\cdots<[n]\qquad\hbox{and}\qquad
1<12<123<1234<\cdots<[n]$$
are both maximal chains
but the first has length $n-1$ while the second has length $n$. Since the posets $L_n$ for $n\geq 4$ are not
graded, they are also not (upper) semimodular.
Indeed, $12$ and $14$ both cover $12\wedge 14 = 1$, but $12\vee 14 = 1234$ does not cover $12$.

For two elements $x\leq y$ in a poset $X$, the {\it interval} $[x,y]$ is the set of all $z\in X$ satisfying
$x\leq z\leq y$. If $X$ has a minimum element $\widehat 0$, then we can define the
the {\it principal (order) ideal generated by $x$}, denoted $\ideal x$, to be the interval $[\widehat 0, x]$.
A poset is said to be {\it locally finite} if every interval is finite.
The {\it M\"obius function} $\mu_X$ of a locally finite poset $X$ is the function from intervals of the poset
to the complex field $\CC$ given by the formulas $\mu_X(x,x) = 1$ for all $x\in X$ and
\newcount\mobiusdef
\mobiusdef=\eqcount
$$\mu_X(x,y) = -\!\!\sum_{x\leq z<y} \mu_X(x,z),\adveq$$
for all $x\leq y$ in $X$, where we have abbreviated $\mu_X\big([x,y]\big)$ by $\mu_X(x,y)$. If the poset $X$
is a lattice, with minimum element $\widehat 0$ and maximum element $\widehat 1$, then
$X = [\widehat 0,\widehat 1]$ and it makes sense to write $\mu_X(X)$ for $\mu(\widehat 0, \widehat 1)$.
In the case that
$X$ is the set of all positive integers, ordered by divisibility, then $\mu_X(m,n) = \mu(n/m)$, where $\mu$
is the classical M\"obius function. Recall that
$\mu(s)=1$ if $s=1$ or $s$ is a product of an even number of distinct
primes, $\mu(s)=-1$ if $s$ is a product of an odd number of distinct primes, and $\mu(s)=0$
if $s$ is divisible by a perfect square. We centre our discussion around the following main result.

\newcount\mobiusformula
\mobiusformula=\thmcount
\proclaim Theorem \advthm. Let $\mu_n = \mu_{L_n}$
be the M\"obius function of the lattice of arithmetic progressions $L_n$.
We have $\mu_0(L_0) = 1$, $\mu_1(L_1) = -1$, and $\mu_n(L_n) = \mu(n-1)$ for $n\geq 2$, where $\mu$ is
the classical M\"obius function.\slug

We now briefly outline the paper. In Section 2, we develop some properties of the number $p_{nk}$
of arithmetic progressions of size $k$ in $[n]$ and show that these quantities arise in a recurrence that
proves
Theorem~{\the\mobiusformula} directly from the definition of the M\"obius function.
In Section 3, we count chains in $L_n$ in order to gain information about the order complex of $L_n$
and derive the same recurrence in a slightly different manner.
We then proceed in Section 4 to study the set of coatoms in $L_n$ in order
to give a general formula for $\mu_n$,
evaluated at an arbitrary interval of $L_n$. As a corollary, we obtain a third proof of
Theorem~{\the\mobiusformula} that is of a rather different nature than the first two proofs.
In Section 5, we explicitly compute the homology groups of the order complex $\Delta_n$ of $L_n$. In Section 6, we
prove that $L_n$ is comodernistic, a property recently introduced by J.~Schweig and R.~Woodroofe that in
particular implies that $\Delta_n$ is shellable for all $n$~\ref{21}.
Lastly, in Section 7, we use lemmas
proved in previous sections to show that $L_n$ is EL-labelable, that $\Delta_n$ is either contractible or has
the homotopy type of a sphere, and that $L_n$ is complemented if and only if $n-1$ is squarefree.

\advsect The number of arithmetic progressions
\hldest{xyz}{}{numprogs}

Our starting point is the number $p_{nk}$ of arithmetic progressions of length $k$ contained in $[n]$.
It was shown in~\ref{11} that for $2\leq k\leq n$,
\newcount\pnk
\pnk=\eqcount
$$p_{nk} = \sum_{r=1}^{\lfloor(n-1)/(k-1)\rfloor} \big(n-(k-1)r\big)
  = n\bigg\lfloor{n-1\over k-1}\bigg\rfloor - {k-1\over 2}
\bigg(\bigg\lfloor{n-1\over k-1}\bigg\rfloor^2 +\bigg\lfloor{n-1\over k-1}\bigg\rfloor \bigg).\adveq$$
(We have halved their formula here, because we consider arithmetic progressions as sets and not as ordered
sequences.)
We also have $p_{n0} = 1$ to count the empty progression as well as $p_{n1}=n$ to count the $n$ singletons.
Values of $p_{nk}$ for small values of $n$ and $k$ are collected in Table 1.
We first derive a formula for the bivariate generating function of $p_{nk}$
(see, e.g.,~\ref{8} for an exhaustive reference on generating functions).

\newcount\genfct
\genfct=\thmcount
\proclaim Lemma \advthm. For integers $n,k\geq 0$, let $p_{nk}$ denote the number of arithmetic progressions
of size $k$ in the interval $[n]$. We have the formula
$$f(z,q) = \sum_{k=0}^\infty \sum_{n=0}^\infty p_{nk} z^n q^k
= {1\over (1-z)^2}\bigg(1-z+zq + \sum_{k=2}^\infty {(zq)^k \over 1-z^{k-1}}\bigg)\adveq$$
for the bivariate generating function of $p_{nk}$.

\goodbreak
\proof The sequences $\big(p_{n0}\big)_{n\geq 0}$ and $\big(p_{n1}\big)_{n\geq 0}$ are $(1,1,1,\ldots)$
and $(0,1,2,\ldots)$ respectively, so that the coefficient of $q^0$ in $f(z,q)$ is $1/(1-z)$ and the
coefficient of $q$ is $z/(1-z)^2$. For $k\geq 2$, there are $n-1$ possible base points and for each base point $a$,
the number of possible step sizes is $\lfloor (n-a)/(k-1)\rfloor$. So
$$\sum_{n=2}^\infty p_{nk}z^n = \sum_{n=2}^\infty \sum_{a=1}^{n-1}\bigg\lfloor {n-a\over k-1} \bigg\rfloor z^n
=\sum_{n=0}^\infty \sum_{a=1}^{n-1}\bigg\lfloor {a\over k-1} \bigg\rfloor z^n
={1\over 1-z}\sum_{n=0}^\infty \bigg\lfloor {n\over k-1} \bigg\rfloor z^n\adveq$$
where we have added the empty terms for $n=0$ and $n=1$ and reversed the order of summation in the second equality.
Note that
$$\eqalign{
\sum_{n=0}^\infty\bigg\lfloor {n\over k-1} \bigg\rfloor z^n&=\sum_{i=1}^\infty\sum_{n=i(k-1)}^\infty z^n\cr
&= \sum_{i=1}^\infty {z^{i(k-1)}\over 1-z}\cr
&= {1\over 1-z}\cdot{1-(1-z^k)\over 1-z^{k-1}}\cr
&= {z^k\over (1-z^{k-1})(1-z)}.\cr
}\adveq$$
Putting everything together, we find that
$$\sum_{k=0}^\infty \sum_{k=0}^\infty p_{nk}z^nq^k
= {1\over 1-z} + {z\over (1-z)^2}q + \sum_{k=2}^\infty {z^k\over (1-z^{k-1})(1-z)^2} q^k,\adveq$$
which simplifies to the formula we were looking for.\slug

\topinsert
\vskip-10pt
$$\vcenter{\vbox{
\eightpoint
\centerline{\smallheader Table 1}
\medskip
\centerline{THE NUMBER $p_{nk}$ OF ARITHMETIC PROGRESSIONS OF SIZE $k$ IN $\{1,2,\ldots,n\}$}
}}$$
\vskip-10pt
$$\vcenter{\vbox{
\eightpoint
\hrule
\medskip
\tabskip=.7em plus.2em minus .5em
\halign{
   $\hfil#$  & $\hfil#$ &  $\hfil#$ & $\hfil#$ & $\hfil#$ & $\hfil#$ & $\hfil#$ & $\hfil#$ & $\hfil#$ &
   $\hfil#$ & $\hfil#$ & $\hfil#$ & $\hfil#$ \cr
   n   & p_{n0} & p_{n1} & p_{n2} & p_{n3} & p_{n4} & p_{n5} & p_{n6} & p_{n7} & p_{n8} & p_{n9} & p_{n(10)}
   & p_{n(11)} \cr
   \noalign{\medskip}
   \noalign{\hrule}
   \noalign{\medskip}
1 & 1 & 1 \cr
2 & 1 & 2 & 1 \cr
3 & 1 & 3 & 3 & 1 \cr
4 & 1 & 4 & 6 & 2 & 1 \cr
5 & 1 & 5 & 10 & 4 & 2 & 1 \cr
6 & 1 & 6 & 15 & 6 & 3 & 2 & 1 \cr
7 & 1 & 7 & 21 & 9 & 5 & 3 & 2 & 1 \cr
8 & 1 & 8 & 28 & 12 & 7 & 4 & 3 & 2 & 1 \cr
9 & 1 & 9 & 36 & 16 & 9 & 6 & 4 & 3 & 2 & 1 \cr
10 & 1 & 10 & 45 & 20 & 12 & 8 & 5 & 4 & 3 & 2 & 1 \cr
11 & 1 & 11 & 55 & 25 & 15 & 10 & 7 & 5 & 4 & 3 & 2 & 1 \cr
 \noalign{\medskip}
 \noalign{\hrule}
    }
}}$$
\vskip-10pt
\endinsert

Because $p_{nk} = 0$ when $k>n$,
the horizontal generating functions $f_n(q)$ are polynomials $\sum_{k=0}^n p_{nk}q^k$. For instance, since $L_0$
through $L_3$ are just boolean lattices (consisting of all subsets of a finite ground set),
we have $f_n(q) = (1+q)^n$. When $n=4$, we have
$f_4(q) = 1 + 4q + 6q^2 + 2q^3 + q^4$, which is irreducible in $\ZZ[q]$ by Cohn's criterion~\ref{5},
since $f_4(10) = 12641$ is prime. It can also be checked computationally that $f_n(q)$ is irreducible
for $5\leq n\leq 10$, and there is no reason to suspect that this polynomial has a neat factorisation for
any larger values of $n$.
As a corollary of the above lemma, we obtain a nice formula for $f_n(1) = |L_n|$, the number of elements in the
lattice.

\newcount\sizeformula
\sizeformula=\eqcount
\proclaim Corollary \advthm. For $n\in \NN$, the poset $L_n$ has
$$|L_n| = 1 + n + \sum_{a=1}^{n-1}\sum_{r=1}^a \tau(r)\adveq$$
elements, where $\tau(r) = \sum_{d\divides r} 1$ is the divisor function.

\proof We write
$$|L_n| = f_n(1) = 1+ n +\sum_{a=1}^{n-1} \sum_{r=1}^a \bigg\lfloor {a\over r}\bigg\rfloor$$
and then apply the elementary identity
$\sum_{k=1}^n\tau(k) = \sum_{k=1}^n \lfloor n/ k\rfloor$.\slug

The sequence $(|L_n|-1)_{n\geq 1}$ appears in the On-line Encyclopedia of Integer Sequences under the
entry A051336.
We now proceed to the first proof of Theorem~{\the\mobiusformula}, which expresses the M\"obius function
of $L_n$ as a recurrence defined in terms of $p_{nk}$.

\medskip\noindent{\it First proof of Theorem~{\the\mobiusformula}.}\enspace
Let $M_n = \mu_n(L_n)$ for short. The case $n=0$ is trivial. For $n\geq 1$,
we must subtract $\mu_k(\emptyset, x)$ for every progression $x\in L_n^* = L_n\setminus\{[n]\}$.
Because $x = \big\{a,a+r,\ldots,a+(k-1)r\big\}$ is a progression, one obtains
an isomorphism of posets between the ideal $\ideal x$ and $L_k$ by relabelling the element $a+ir$ with
$i+1$ for $0\leq i<k$. Hence $\mu_n(\emptyset, x) = \mu_k(\emptyset, [k]) = M_k$, and since there are
$p_{nk}$ progressions of size $k$ in $L_n$, we have the recurrence
\newcount\firstthmrecurrence
\firstthmrecurrence=\eqcount
$$M_n = -\sum_{x\in L_n^*} \mu_n(\emptyset, x) = -\sum_{k=0}^{n-1} M_k p_{nk}.\adveq$$
We can then compute $M_1 = -1$ and $M_2 = 1 = \mu(1)$. For $n>2$ we now proceed by strong induction;
suppose that $M_k = \mu(k-1)$ for all $2\leq k<n$. We expand the above recurrence to
$$\eqalign{
M_n &= - \bigg(M_0 p_{n0} + M_1 p_{n1} + \sum_{k=2}^{n-1} M_k p_{nk}\bigg)\cr
&= -\bigg(1 - n + \sum_{k=2}^{n-1} \mu(k-1) \sum_{r=1}^{\lfloor(n-1)/(k-1)\rfloor} \big(n-(k-1)r\big)\bigg)\cr
&= -\bigg(1 - n -\mu(n-1) + \sum_{k=1}^{n-1} \mu(k) \sum_{r=1}^{\lfloor(n-1)/k\rfloor} \big(n-kr\big)
\bigg)\cr
}\adveq$$
and sum over all possible values of $kr$ by setting $m=kr$ and summing over divisors $d$ of $m$, for
$1\leq m\leq n-1$. This gives
$$M_n = -\bigg(1 - n -\mu(n-1)+ \sum_{m=1}^{n-1} \sum_{d\divides m} \mu(d) (n-m)\bigg).\adveq$$
But $\sum_{d\divides m}\mu(d) = 0$ when $m>1$ and when $m=1$, the summation equals $n-1$. After cancellation,
we see that the right-hand side equals $\mu(n-1)$, which is what we wanted to show.\slug

\advsect Chains and the order complex
\hldest{xyz}{}{chains}

An {\it abstract simplicial complex} is a set system $\Delta$
on a vertex set $V$ containing every singleton subset of $V$ and
with the property that for every set $F\in \Delta$, all subsets of $F$ also belong to $\Delta$. The
elements of $\Delta$ are called {\it faces}, and the {\it dimension} of a face $F$ is defined to be $|F|-1$.
A face is said to be {\it maximal} if it is not strictly contained in another face, and the dimension of $\Delta$
is the maximum dimension of a (maximal) face in $\Delta$.
For our purposes, simplicial complexes will contain the empty set, a face of dimension $-1$. We will
require various notions from topology in this section. Any definitions that we do not recall
here can be found in any introductory textbook, such
as~\ref{19}, for example.

A {\it chain} of length $k$ in a poset $X$ is a set $\{x_1, x_2, \ldots, x_{k+1}\}\subseteq X$ such that
$x_1 < x_2 <\cdots <x_{k+1}$; so a chain of length $0$ is a singleton set.
One can associate a simplicial complex, called the {\it order complex}, to any lattice (with bottom element
$\hat 0$ and top element $\hat 1$)
by taking $L\setminus\{\hat 0,\hat 1\}$ as the vertex set and letting the faces be chains in this modified poset.
Let $L_n'$ denote the poset $L_n$ with the minimum element $\emptyset$ as well as the maximum element $[n]$
removed. Note that chains in $L_n'$ of length $k-2$ are in bijection with chains of length $k$ in $L_n'$
that contain both $\emptyset$ and $[n]$, which we shall count in the next lemma.

\newcount\chainrecurrence
\chainrecurrence=\eqcount
\newcount\chainlemma
\chainlemma=\thmcount
\proclaim Lemma \advthm. The number $b_{nk}$ of chains of length $k$ in $L_n$ that contain $\emptyset$
and $[n]$ satisfies the recurrence
$$b_{nk} = \sum_{i=1}^{n-1} p_{ni} b_{i(k-1)},\adveq$$
for $2\leq k\leq n$, with $b_{n1} = 1$ for all $n$ and $b_{nk} = 0$ whenever $k>n$.

\goodbreak
\proof
The case $k=1$ is trivial and it is clear that $b_{nk}$ should be zero for $k>n$. In the other
cases, we are counting chains $\emptyset \subset x_1\subset \cdots x_k$ (we require
strict inclusion here). We split up the cases
by the second-greatest element $x_k$ of the chain. It is clear that the subchain
$\{\emptyset,x_1,\ldots,x_k\}$ is a chain containing both the maximum and minimum element of
the ideal $\ideal x_k$, which is isomorphic to $L_i$, where
$i$ is the size of $x_k$ (as a set). Thus the number of such chains is $b_{m(k-1)}$. There were $p_{ni}$
choices for the element of size $i$, and summing over all possible $i$ gives the recurrence above.\slug

\topinsert
\vskip-10pt
$$\vcenter{\vbox{
\eightpoint
\centerline{\smallheader Table 2}
\medskip
\centerline{THE NUMBER $b_{nk}$ OF CHAINS OF LENGTH $k$ IN $L_n'$}
}}$$
\vskip-10pt
$$\vcenter{\vbox{
\eightpoint
\hrule
\medskip
\tabskip=.7em plus.2em minus .5em
\halign{
   $\hfil#$  &  $\hfil#$ & $\hfil#$ & $\hfil#$ & $\hfil#$ & $\hfil#$ & $\hfil#$ & $\hfil#$ &
   $\hfil#$ & $\hfil#$ & $\hfil#$ & $\hfil#$ & $\hfil#$ \cr
   n   & b_{n1} & b_{n2} & b_{n3} & b_{n4} & b_{n5} & b_{n6} & b_{n7} & b_{n8} & b_{n9} & b_{n(10)}
   & b_{n(11)} \cr
   \noalign{\medskip}
   \noalign{\hrule}
   \noalign{\medskip}
1 &  1 \cr
2 & 1 & 2 \cr
3 & 1 & 6 & 6 \cr
4 & 1 & 12 & 24 & 12 \cr
5 & 1 & 21 & 68 & 72 & 24 \cr
6 & 1 & 32 & 144 & 244 & 180 & 48 \cr
7 & 1 & 47 & 283 & 666 & 764 & 432 & 96 \cr
8 & 1 & 64 & 486 & 1510 & 2436 & 2164 & 1008 & 192 \cr
9 & 1 & 85 & 799 & 3117 & 6534 & 8028 & 5816 & 2304 & 384 \cr
10 & 1 & 109 & 1232 & 5860 & 15368 & 24524 & 24516 & 15040 & 5184 & 768 \cr
11 & 1 & 137 & 1838 & 10418 & 33049 & 65402 & 84284 & 70992 & 37760 & 11520 & 1536 \cr
 \noalign{\medskip}
 \noalign{\hrule}
    }
}}$$
\vskip-10pt
\endinsert

For small values of $n$ and $k$, the values $b_{nk}$ are displayed in Table 2.
Note that if in the recurrence~\refeq{\the\chainrecurrence} we replace $p_{nk}$ with ${n\choose k}$, we
obtain the array of numbers $k!\big\{\!{n\atop k}\!\big\}$, where
$\big\{\!{n\atop k}\!\big\}$ is a Stirling number
of the second kind (see, e.g.,~\ref{13}).
These numbers count the number of ways to partition $n$ numbers into $k$ nonempty subsets,
and for each such partition $S_1, S_2, \ldots, S_n$, we obtain $k!$ chains in the boolean lattice that
contain both $\emptyset$ and $[n]$ (for each
permuation $\sigma$ in $\frak S_n$, we have the chain $\{\emptyset, S_{\sigma(1)},
S_{\sigma(1)} \cup S_{\sigma(2)}, \ldots, [n]\}$).

Returning to our numbers $b_{nk}$, we see that for $-1\leq k\leq n-2$, the number of $k$-dimensional
faces of $\Delta_n$ is $b_{n(k+2)}$. Hence $\Delta_n$ is an $(n-2)$-dimensional simplicial complex. Let
$\tilde\chi(\Delta_n) = \chi(\Delta_n)-1$ be the reduced Euler characteristic of the order complex. We have
$$\tilde\chi(\Delta_n) = \sum_{k=1}^n (-1)^k b_{nk}\adveq$$
for $n\geq 1$.
Using the fact that the M\"obius function of a poset with a maximum and minimum element artificially adjoined
equals the reduced Euler characteristic of its order
complex, we obtain an alternative proof of Theorem~{\the\mobiusformula}.

\medskip\noindent{\it Second proof of Theorem~{\the\mobiusformula}.}\enspace Let
$M_n = \sum_{k=1}^n (-1)^k b_{nk}$.
We compute $M_0 = 1$ and $M_1 = -1$ by hand. To complete the proof, it suffices to show that
$\tilde\chi(\Delta_n)= M_n = \mu(n-1)$
for all $n\geq 2$. The base case $M_2 = 1$ follows from a direct computation, and for $n>2$, we have
$$\eqalign{
M_n &= \sum_{k=1}^n (-1)^k b_{nk}\cr
&= -1 + \sum_{k=2}^n (-1)^k \sum_{i=1}^{n-1} p_{ni}b_{i(k-1)}\cr
&= -1 + \sum_{i=1}^{n-1} p_{ni} \sum_{k=2}^n (-1)^k b_{i(k-1)},
}\adveq$$
by Lemma~{\the\chainlemma}. We can pull out one of the $-1$ factors and reindex to obtain
$$M_n = -\bigg(1+ \sum_{i=1}^{n-1} p_{ni} \sum_{k=1}^{i}(-1)^k b_{ik}
\bigg).\adveq$$
Note that the upper index in the inner summation has been changed to
$i$, since $b_{ik}=0$ when $k>i$. By the induction
hypothesis, this inner sum is $M_i$, so
$$M_n = -\bigg(1 + \sum_{i=1}^{n-1} p_{ni} M_i\bigg) = \sum_{i=0}^{n-1} p_{ni}M_i,\adveq$$
which is the recurrence~\refeq{\the\firstthmrecurrence} we encountered in the first proof
of this theorem. The rest of the proof proceeds
exactly as before.
\slug

\advsect Coatoms
\hldest{xyz}{}{coatoms}

We now set out to compute $\mu_n(x_1, x_2)$ for arbitrary progressions $x_1$ and $x_2$ in $L_n$. Towards
this goal, we will need to study the {\it coatoms} of $L_n$, the elements covered by $[n]$. It turns
out that we can give an explicit description of the set of coatoms in $L_n$. In the following lemma, we
use the notation $j\divides k$ to indicate that $k$ is an integer multiple of $j$.

\newcount\coatoms
\coatoms=\thmcount
\proclaim Lemma \advthm. Let $A_n\subseteq L_n$ be the set of coatoms. We have $A_1 = \{\emptyset\}$,
$A_2 = \{1, 2\}$, and $A_3 = \{12, 13, 23\}$. For $n\geq 4$, we have $A_n = B_n \cup C_n$, where
$B_n = \{12\cdots(n-1), 23\cdots n\}$, and
$$C_n = \cases{\{1n\}, & if $n-1$ is prime;\cr
\big\{\{1, 1+p, 1+2p, \ldots, n\} : p\ \hbox{\rm prime},\ p\divides n-1\big\},& otherwise.\cr}\adveq$$
In particular, the size of $A_n$ is $\omega(n-1) +2$, where $\omega(n)$ is the number of distinct prime
divisors of $n$.

\proof The small cases are easily computed explicitly.
When $n\geq 4$ there are only two elements of size $n-1$, and the fact that they are coatoms is
obvious. Now any element that does not contain both $1$ and $n$ cannot be a coatom, since an element of
$B_n$ would contain it. The progressions that contain $1n$ are of the form $x_d = \{1, d+1, 2d+1, \ldots, n\}$
for divisors $d$ of $n-1$, but note that if $d$ is composite, then $x_d$ is contained in $x_{d'}$ for any
$d'$ dividing $d$. Hence the remaining coatoms are the progressions with prime steps, implying that
$C_n$ is of one of the two forms above.\slug

Every element in $L_n$ is contained in some coatom, but not all elements can be expressed as a {\it meet}
of coatoms. The next lemma shows that in $L_n$, if an element can be expressed as a meet of coatoms,
then this representation is unique.

\newcount\uniquemeet
\uniquemeet=\thmcount
\proclaim Lemma \advthm. Let $L_n$ be the lattice of arithmetic progressions and let $A_n\subseteq L_n$ be
the set of coatoms. If $x\in L_n$ can be expressed as $x= \bigwedge_{s\in S} s$ for some $S\subseteq A_n$,
then $S$ is uniquely determined by $x$.

\proof If $x = \emptyset$, the only possibility is to take $S = A_n$, since omitting
one of $12\cdots(n-1)$ or $23\cdots n$ would cause one of the elements $1$ or $n$ to appear in the meet,
and omitting the progression with base point $1$, step size $p$ (a prime dividing $n-1$),
and end point $n$ will cause the $p-1$ elements
$$1 + {n-1\over p}, 1+ {2(n-1)\over p}, \ldots,  n-{n-1\over p} $$
to appear in the meet.

Now suppose that $x$ is nonempty and we can write out the elements of
$x = \big\{a, a+r, \ldots, a+(k-1)r\big\}$. We will
consider the possible step sizes $r$. When $r=1$, $x$ is either $12\cdots(n-1)$, $23\cdots n$, or
$23\cdots(n-1)$ and in all three cases it is clear that there is only one representation of $x$ as the meet
of coatoms. For $r>1$, we find that $r$ must be the least common multiple of some primes dividing $n-1$,
and there is only one way to express $r$ as a least common multiple of distinct primes, thus uniquely determining
the coatoms with prime step size that are in $S$. Lastly, note that $12\cdots (n-1)$ is in $S$ if and only
if $a = 1$ and $23\cdots n$ is in $S$ if and only if $a+(k-1)r = n$.\slug

These properties of the set of coatoms in $L_n$ implies a general formula for computing $\mu_n(x,[n])$.

\newcount\generalthm
\generalthm=\thmcount
\proclaim Theorem \advthm. Let $x$ be an arbitrary progression in $L_n$. For all $[n]\neq x \in L_n$,
$$\mu_n(x, [n]) = \cases{(-1)^k, & if $x$ is the meet of $k$ coatoms;\cr 0,& if $x$ is not a meet
of coatoms.\cr}\adveq$$

\proof Note that $x$ is the minimum element of the interval $L = \big[x, [n]\big]$; the subset $S\subseteq L_n$
of coatoms whose meet equals $x$ is contained in this interval. By the cross-cut theorem~\ref{20},
$$\mu_n(x, [n]) = \sum_{k=1}^{|A_n|} (-1)^k N_k,\adveq$$
where $N_k$ is the number of subsets of $S$ of size $k$ whose meet is $S$. By Lemma~{\the\uniquemeet},
$N_{|S|} = 1$ and $N_k = 0$ for all $k\neq |S|$, proving the theorem.\slug

It is easy to tell if a given progression $x$ is a meet of coatoms, since such $x$ have a very specific form.
In particular, $x$ is a meet of coatoms of $L_n$ if and only if
$$x \cap \{2,\ldots, n-1\} = (1+d\ZZ) \cap \{2,\ldots, n-1\}$$
for some divisor $d$ of $n-1$.
One can then work out the number of elements in the meet representation by taking the prime decomposition of $d$
and checking whether $1$ or $n$ (or both or neither) are included in $x$.
Let $\omega(n)$ be the number of distinct primes dividing an integer $n$ and let $S$ denote the set of
progressions $x$ with $\mu_n(x,[n])\neq 0$. Lemma~{\the\coatoms} and Theorem~{\the\generalthm}
together imply that there are exactly
$2^{\omega(n-1)+2}$ such elements $x$ in $L_n$. Since every squarefree divisor of $n-1$ contributes exactly
four progressions to the set $S$, we can prove the elementary identity
$\sum_{d\divides n}\big|\mu(d)\big| = 2^{\omega(n)}$ by counting $S$ in two ways.

Since the ideal $\ideal x\subseteq L_n$ is isomorphic to $L_m$ for any progression $x$ of size $m$,
Theorem~{\the\generalthm}
immediately implies a general method for computing the M\"obius function of an arbitrary interval.

\proclaim Corollary \advthm. Let $x_1$ and $x_2$ be arbitrary elements of $L_n$ and let $C$ be the set of elements
covered by $x_2$. We have
$$\mu_n(x_1, x_2) = \cases{(-1)^k, & if $x_1$ is the meet of $k$ elements of $C$;\cr 0,& if $x_1$ is not a meet
of elements of $C$.\noskipslug\cr}\adveq$$

This corollary tells us that the M\"obius function of $L_n$ takes values in $\{0,\pm 1\}$ no matter the
interval at which it is evaluated. Posets with this
property are sometimes called {\it totally unimodular} (see, e.g.,~\ref{15}).
Theorem~{\the\generalthm} also allows us to give a third proof of Theorem~{\the\mobiusformula}.

\medskip\noindent{\it Third proof of Theorem~{\the\mobiusformula}.}\enspace We take $n\geq 4$; smaller cases
can easily be worked out explicitly.
First suppose that $n-1$ is squarefree, equalling the product of distinct primes $p_1,p_2,\ldots, p_k$,
so that $\mu(n-1) = (-1)^k$.
The claim is that for any nonempty progression $x\in L_n$, there is some coatom that does not contain $x$.
If $x$ contains either $1$ or $n$, then one of the two progressions in $L_n$ of size $n-1$ does not contain $x$.
Otherwise, $x$ contains some integer $1+m$ for $1\leq m \leq n-2$.
Since $m< n-1 = \lcm(p_1,p_2\ldots, p_k)$, there is some
prime $p_i$ that does not divide $m$, hence $1+m$ is not contained in the coatom of step size $p_i$.
There are $k+2$ coatoms in $L_n$, so
Theorem~{\the\generalthm} can be applied to give $\mu_n(L_n) = (-1)^{k+2} = (-1)^k = \mu(n-1)$.

Now assume that $n-1$ is divisible by $p^2$ for some prime $p$. Since the integer $(n-1)/p$ is divisible by
every prime dividing $n-1$, the element $1+(n-1)/p$ belongs to every coatom of $L_n$. So $\emptyset$ cannot
be expressed as a meet of coatoms and $\mu_n(L_n) = 0$.\slug

\advsect Homology groups of the order complex
\hldest{xyz}{}{homology}

Although less direct than the first two proofs we supplied, the
proof of Theorem~{\the\mobiusformula} given in the previous section reveals much of the internal
structure of $L_n$. We now show that it can be reinterpreted to give a complete characterisation of the homology
groups of $\Delta_n$, a strictly stronger result than Theorem~{\the\mobiusformula}.
A simplicial complex $\Delta$, as we have defined it, is simply a set system,
but $\Delta$ can be embedded in Euclidean space to give rise to a topological space
$|\Delta|$ called its {\it geometric realisation}. We will sometimes abuse notation and ascribe
topological properties of $|\Delta|$ to $\Delta$.
The reduced Euler characteristic of an $n$-dimensional
simplicial complex $\Delta$ can also be expressed as the alternating sum
$$\tilde\chi(\Delta) = \tilde\chi(|\Delta|) = \sum_{i=0}^n (-1)^i \rank \tilde H_i(|\Delta|, \ZZ),\adveq$$
where $\tilde H_i(|\Delta|,\ZZ)$ is the $i$th reduced homology group of the topological space $|\Delta|$
(whenever we refer to a homology group, we shall understand reduced homology group).

To derive the homology groups of $L_n$, we will require the notion of cross-cuts. A {\it cross-cut} $C$ of a
lattice $L$ (with maximum $\widehat 1$ and minimum $\widehat 0$)
is a subset of $L$ not containing either of $\widehat 1$ and $\widehat 0$ such that no two
elements of $C$ are comparable and every maximal chain in the lattice contains some element of $C$. A subset
$S$ of $L$ is said to be {\it spanning} if the join of all its elements is $\widehat 1$ and the
meet of all its elements is $\widehat 0$. For a cross-cut $C$ of a lattice $L$, we can define a simplicial
complex $\Delta(C)$ whose vertices are the elements of $C$ and whose faces are given by subsets of $C$
that are {\it not} spanning. A paper of J.~Folkman~\ref{9}
showed that $\tilde H_i\big(\Delta(C), \ZZ) \cong
\tilde H_i\big(\Delta, \ZZ)$ for all $i$, where $\Delta$ is the order complex of $L$. We use this to derive
the homology groups of $\Delta_n$.

\newcount\homology
\homology=\thmcount
\proclaim Lemma \advthm. For $n\geq 4$, let $L_n$ be the lattice of arithmetic progressions and let $\Delta_n$
be the order complex of $L_n' = L_n\setminus\{\emptyset, [n]\}$.
Let $\tilde H_i(\Delta_n, \ZZ)$ be the $i$th reduced homology group of $\Delta_n$.
If $n-1$ is squarefree and equal to the product of $k$ distinct primes, then
$$\tilde H_i(\Delta_n, \ZZ) = \cases{\ZZ, & if $i=k$; \cr 0, & otherwise.\cr}.$$
If $n-1$ is not squarefree, then all the homology groups of $\Delta_n$ are trivial.

\proof Let $C$ be the set of coatoms of $L_n$, whose explicit construction was given by Lemma~{\the\coatoms}.
Let $k = \omega(n-1) = k$, so that $|C| = k+2$. If $n-1$ is squarefree, then as we saw earlier in the third
proof of Theorem~{\the\mobiusformula}, we can express $\emptyset$ as a meet of elements of $C$, so $C$ is a
spanning set. However, any proper subset $C'$ of $C$ is not spanning, since if $c_i$ is the element of $C$ that
is not in $C'$, then we can build a chain $\emptyset \subset \cdots \subset c_i\subset [n]$ that does not
contain an element of $C'$. So every subset of $C$ with cardinality $k+1$ is an element of the abstract simplicial
complex $\Delta(C)$, i.e., $\Delta(C)$ is the boundary of a $(k+1)$-dimensional simplex, whose $k$th homology
group is $\ZZ$ and whose other reduced homology groups are all trivial.

When $n-1$ is not squarefree, the construction we gave in the third proof of Theorem~{\the\mobiusformula} shows
that $\emptyset$ is not the meet of the elements of $C$, which means that $C$ itself does not span. Hence
$\Delta(C)$ is the $(k+1)$-dimensional simplex, including its interior, all of whose
reduced homology groups are trivial.\slug

We will use Lemma~{\the\homology} later on to prove the stronger fact that $\Delta_n$ has the homotopy type
of a sphere when $n-1$ is squarefree.

\advsect Left-modularity and comodernism
\hldest{xyz}{}{comodernism}

An element $m$ in a lattice $L$ is {\it left-modular in $L$}
if for all $x<y \in L$, $(x\vee m) \wedge y = x \vee (m\wedge y)$.
A lattice $L$ is {\it comodernistic} if every interval $[x,y]\subseteq L$ has a coatom which is left-modular in
$[x,y]$. The aim of this section is to show that $L_n$ is comodernistic.
To do so, we will make use of two of the lemmas in the paper of J.~Schweig
and R.~Woodroofe that introduced the definition of comodernism.

\parenproclaim Lemma A ({\rm\ref{21}}, Lemma 2.12).
Let $m$ be a coatom of the lattice $L$. Then $m$ is left-modular in $L$ if and only if for every $y \in L$ with
$y \not\leq m$, $y$ covers $m \wedge y$.\slug

\parenproclaim Lemma B ({\rm\ref{21}}, Lemma 4.1).
Let $L'$ be a sublattice of a lattice $L$. If $m\in L'$ is a left-modular coatom in $L$, then $m$ is also left-modular
in $L'$.\slug

Note that we have modified these lemmas slightly to suit our notation and use case; in particular, the original
version of Lemma B requires only that $L'$ be a meet subsemilattice. We begin with a small lemma.

\newcount\bnatoms
\bnatoms=\thmcount
\proclaim Lemma \advthm. For $n\ge 1$, the elements $12\cdots(n-1)$ and $23\cdots n$ are left-modular in $L_n$.

\proof
Without loss of generality, let $m = 12 \cdots(n-1)$; the case where $m=23\cdots n$ is symmetric.
Let $y \in L_n$ be such that $y \not\leq m$, so it must be that $n \in y$,
hence $m\wedge y = y \setminus \{n\}$ which is covered by $y$. By Lemma A, this shows that $m$ is left-modular.\slug

We are now able to show that $L_n$ is comodernistic for all $n$. For brevity of notation, in the following proof
we let $\filter_k x$ denote the {\it principal filter} of
$x\in L_k$; that is, $\filter_k x = \{y\in L_k : x \leq y\}$.

\proclaim Theorem \advthm. For all $n\geq 0$, the lattice $L_n$ is comodernistic.

\proof Let $[x,y]$ be an interval in $L_n$. We once again employ the fact that
$\ideal y$ is isomorphic to $L_k$ where $k = |y|$.
This isomorphism sends $[x,y]$ to the interval $\filter_k x \subseteq L_k$, so it suffices to show that,
for all $k \ge 1$ and $x\in L_k$, the principal filter $\filter_k x$ contains a coatom which is left-modular (in
the filter). Let $A_k = B_k \cup C_k$ be the coatoms of $L_k$, with $B_k$ and $C_k$
defined as in Theorem~{\the\coatoms}.
Clearly, the coatoms of $\filter_k x$ are a subset of $A_k$.
If $\filter_k x\cap B_k \neq \emptyset$,
then by Lemma~{\the\bnatoms}, $\filter_k x$ contains a coatom which is left-modular
in all of $L_k$, and by Lemma B it is also a left-modular coatom in $\filter_k x$.
If ${\filter_k x}\cap B_k$ is empty, then the progression $x$ must contain both $1$ and $k$, so
$\filter_k x\subseteq \filter_k 1k$ and in particular, every coatom of $\filter_k x$ is also a coatom
of $\filter_k 1k$. By another application of Lemma B, we may reduce our proof to showing that every coatom
of $\filter_k 1k$ is left-modular in this filter.

The coatoms of $\filter_k 1k$ are precisely the elements in $C_k$.
If $k-1$ is prime, Lemma~{\the\coatoms} tells us that $1k$ is a coatom, so $\filter_k 1k$ contains only the two
elements $1k$ and $[k]$, the former of which is trivially left-modular in this interval.
On the other hand, let $k-1$ be composite and let $m$ be a coatom of $\filter_k 1k$; by Lemma~{\the\coatoms}, $m$
is of the form $\{1, 1+p, ... , k\}$ for some $p$ dividing $k-1$.
If $y\in \filter_k 1k$ satisfies $y \not\leq m$, then $y = \{1, 1+r, ... , k\}$
where $r$ divides $k-1$ and $p$ does not divide $r$.
So $m\wedge y = \{1, 1+s, ... , k\}$ where $s = \lcm(r, p) = rp$,
hence $m\wedge y$ is covered by $y$ and we conclude that $m$ is left-modular by Lemma A.\slug

\advsect EL-labelability, homotopy type, and complements
\hldest{xyz}{}{homotopy}

We now use the lemmas of the previous sections to demonstrate further properties of $L_n$. Here we show
that $L_n$ is EL-labelable, that $\Delta_n$ is either homotopy equivalent to a point or a sphere, and that $L_n$
is complemented if and only if $n-1$ is squarefree.

\medskip\boldlabel EL-labelability.
Given a lattice $L$, let $E(L)$ be the set of all $(x,y)\in L$ such that $y$ covers $x$; thus $E(L)$ is the edge set of
the Hasse diagram of $L$. We say that a function $\lambda:E(L)\to \ZZ$ is an {\it ER-labeling} (or {\it edge-rising
labeling)} if for every interval
$[x,y]\subseteq L$, there is a unique maximal chain $x = x_0 < x_1 < \cdots < x_s = y$ with increasing
labels, that is, with
$$\lambda(x_0, x_1) < \lambda(x_1,x_2) < \cdots < \lambda(x_{s-1},x_s).$$
Let $\ZZ^*$ denote the set of all finite sequences of integers.
One defines a lexicographic partial order $\preceq$
on $\ZZ^*$ by declaring $(a_1, \ldots, a_m) \preceq (b_1, \ldots, b_n)$
if either $a_i = b_i$ for $1\leq i\leq m$ and $m\leq n$ or else $a_i<b_i$ for the smallest $i$ with $a_i\ne b_i$.
Note that the function $\lambda$ defines a map $\bar\lambda$ from chains in $L$ to tuples of positive integers;
namely if $c$ is the chain formed by $x_0 < x_1 < \cdots < x_s$, then
$$\bar\lambda(c) = \bigl(\lambda(x_0, x_1) , \lambda(x_1,x_2), \ldots ,\lambda(x_{s-1},x_s)\bigr).$$
Let $\lambda$ be an ER-labeling with the further property that for all $[x,y]$, the unique increasing maximal chain $m$
has $\bar\lambda(m) \preceq \bar\lambda(m')$ for all other maximal chains $m'$ in $[x,y]$. Such an ER-labeling is
called an {\it EL-labeling} (or {\it edge-lexicographic} labeling). A lattice that admits an ER-labeling is
said to be {\it ER-labelable} and one that admits an EL-labeling is {\it EL-labelable}.

A paper of T.~Li showed that comodernistic lattices are EL-labelable~\ref{18}, so in particular we find
that for all $n\geq 0$, $L_n$ is EL-labelable. The looser property of ER-labelability is useful in certain enumerative
problems. For example, it has been shown that the zeta and M\"obius transforms
for ER-labelable posets $P$ can be computed in at most $|E(P)|$
elementary arithmetic operations~\ref{17}.

\medskip\boldlabel Homotopy type. A simplicial complex $\Delta$ is
{\it nonpure shellable} if its maximal faces can be given an order $C_1, C_2, \ldots, C_m$ such that for all
$2\leq k\leq m$, the maximal faces in the complex $\bigl(\bigcup_{i=1}^{k-1} C_i\bigr) \cap C_k$
all have dimension $\dim C_k -1$.
The earliest treatment of nonpure shellable complexes was
carried out by A.~Bj\"orner and M.~L.~Wachs in~\ref{2} and~\ref{3}; Corollary 13.3 of the latter
asserts that a nonpure shellable complex is homotopy equivalent to a wedge of spheres. Proposition 2.3 of an earlier
paper by the same authors~\ref{1} states that EL-labelable posets are nonpure
shellable, so $\Delta_n$ is homotopy
equivalent to a wedge of spheres. In fact, $\Delta_n$ is either contractible or homotopy equivalent to a single sphere,
as the following strengthening of Lemma~{\the\homology} shows.

\proclaim Theorem \advthm. Let $\Delta_n$ be the order complex of the lattice of arithmetic progressions $L_n$.
If $n-1$ is not squarefree, then $\Delta_n$ is contractible. Otherwise, $\Delta_n$ has the homotopy type of
$S^k$, where $k$ is the number of distinct primes dividing $n-1$.

\proof We already know, from the above discussion, that $\Delta_n$ is homotopy equivalent to a wedge of spheres.
If the wedge product consisted of more than one sphere, then the sum over the ranks of the reduced homology groups of
$\Delta_n$ would be greater than $1$. But by Lemma~{\the\homology}, this sum equals $0$ when $n-1$ is not squarefree,
in which case $\Delta_n$ must have the homotopy type of a point, and when $n-1$ is squarefree it equals $1$,
meaning that there exactly one sphere in the wedge product.\slug

\medskip\boldlabel Complements. We finish with a miscellaneous result about complements in $L_n$.
A lattice $L$ with maximum element $\widehat 1$ and minimum element
$\widehat 0$ is said to be {\it complemented} if for all $x\in L$,
there exists $y\in L$ such that $x\vee y = \widehat 1$ and $x\wedge y = \widehat 0$. The elements
$x$ and $y$ are called
{\it complements} of one another, and if we remove the condition that $x\wedge y = \widehat 0$,
then $x$ and $y$ are said to be {\it upper semicomplements}. The next theorem gives
a necessary and sufficient condition for $L_n$ to be complemented.

\newcount\complements
\complements=\thmcount
\proclaim Theorem \advthm. Let $n\geq 2$. The lattice $L_n$ is complemented if and only if $n-1$ is squarefree.
In particular, if $n-1$ is not squarefree, there exists an element $x\notin\{\emptyset, [n]\}$
of $L_n$ whose only upper semicomplement
is $[n]$.

\proof For the ``if'' direction, we note that when $n-1$ is squarefree, we have
$\mu_n(L_n)\neq 0$, which, by a theorem of H.~H.~Crapo~\ref{6}, implies that $L_n$ is
complemented. For the converse, suppose that $n-1$ is divisible by
$p^2$ for some prime~$p$. Consider the progression
$$x = \bigg\{ 1+{n-1\over p}, 1+{2(n-1)\over p}, \ldots, n-{n-1\over p}\bigg\},$$
which has length $p-1$ and is thus not empty.
Note that any $x'\in L$ satisfying $x'\vee x = [n]$ must contain both $1$ and $n$ and the step size $r$
must be coprime to $(n-1)/p$. We also know that $r$
must divide $n-1$. But the only such integer $r$ is $1$, in which case we see that $x$ must be $[n]$.\slug

\section Acknowledgements
\hldest{xyz}{}{acks}

We would like to thank Andrew Granville for helpful conversations
and for giving us the key idea in solving the recurrence that appears in the first
two proofs of Theorem~{\the\mobiusformula}.
We thank
Daniel Wise for enjoyable conversations and for his advice on the topological section in particular.
We are grateful also to Jukka Kohonen for answering a question of the first author on MathOverflow and for
the enlightening discussions that followed via email.
Lastly, we thank our homies Amanda Gu and Rosie Zhao for technical feedback on preliminary drafts.

\section References
\hldest{xyz}{}{refs}

\parskip=0pt
\hyphenpenalty=-1000 \pretolerance=-1 \tolerance=1000
\doublehyphendemerits=-100000 \finalhyphendemerits=-100000
\frenchspacing
\def\bref#1{[#1]}
\def\beginref{\noindent
}
\def\endref{\medskip}
\vskip\parskip

\beginref \parindent=20pt\item{\bref{1}}
\hldest{xyz}{}{bib1}%
Anders Bj\"orner
and Michelle Lynn Wachs,
``On lexicographically shellable posets,''
{\sl Transactions of the American Mathematical Society}\/
{\bf 277}
(1983),
323--341.
\endref

\beginref \parindent=20pt\item{\bref{2}}
\hldest{xyz}{}{bib2}%
Anders Bj\"orner
and Michelle Lynn Wachs,
``Shellable nonpure complexes and posets. I,''
{\sl Transactions of the American Mathematical Society}\/
{\bf 348}
(1996),
1299--1327.
\endref

\beginref \parindent=20pt\item{\bref{3}}
\hldest{xyz}{}{bib3}%
Anders Bj\"orner
and Michelle Lynn Wachs,
``Shellable nonpure complexes and posets. II,''
{\sl Transactions of the American Mathematical Society}\/
{\bf 349}
(1997),
3945--3975.
\endref

\beginref \parindent=20pt\item{\bref{4}}
\hldest{xyz}{}{bib4}%
Anders Bj\"orner,
``A cell complex in number theory,''
{\sl Advances in Applied Mathematics}\/
{\bf 46}
(2011),
71--85.
\endref

\beginref \parindent=20pt\item{\bref{5}}
\hldest{xyz}{}{bib5}%
John Brillhart,
Michael Filaseta,
and Andrew Odlyzko,
``On an irreducibility theorem of A.~Cohn,''
{\sl Canadian Journal of Mathematics}\/
{\bf 33}
1055--1059.
\endref

\beginref \parindent=20pt\item{\bref{6}}
\hldest{xyz}{}{bib6}%
Henry Howland Crapo,
``M\"obius inversion in lattices,''
{\sl Archiv der Mathematik}\/
{\bf 19}
(1969),
595--607.
\endref

\beginref \parindent=20pt\item{\bref{7}}
\hldest{xyz}{}{bib7}%
Richard Ehrenborg,
Likith Govindaiah,
Peter Seho Park,
and Margaret Readdy,
``The van der Waerden complex,''
{\sl Journal of Number Theory}\/
{\bf 172}
(2017),
287--300.
\endref

\beginref \parindent=20pt\item{\bref{8}}
\hldest{xyz}{}{bib8}%
Philippe Flajolet
and Robert Sedgewick,
{\sl Analytic Combinatorics}
(New York:
Cambridge University Press,
2009).
\endref

\beginref \parindent=20pt\item{\bref{9}}
\hldest{xyz}{}{bib9}%
Jon Folkman,
``The homology groups of a lattice,''
{\sl Journal of Mathematics and Mechanics}\/
{\bf 15}
(1966),
631--636.
\endref

\beginref \parindent=20pt\item{\bref{10}}
\hldest{xyz}{}{bib10}%
Harry Furstenberg,
``On the infinitude of primes,''
{\sl American Mathematical Monthly}\/
{\bf 62}
(1955),
353.
\endref

\beginref \parindent=20pt\item{\bref{11}}
\hldest{xyz}{}{bib11}%
Marcel Kieren Goh
and Rosie Y Zhao,
``Arithmetic subsequences in a random ordering of an additive set,''
{\sl arXiv preprint 2012.12339}\/
(2020).
\endref

\beginref \parindent=20pt\item{\bref{12}}
\hldest{xyz}{}{bib12}%
Solomon Wolf Golomb,
``A connected topology for the integers,''
{\sl The American Mathematical Monthly}\/
{\bf 66}
(1959),
663--665.
\endref

\beginref \parindent=20pt\item{\bref{13}}
\hldest{xyz}{}{bib13}%
Ronald Lewis Graham,
Donald Ervin Knuth,
and Oren Patashnik,
{\sl Concrete Mathematics}
(Reading:
Addison--Wesley,
1989).
\endref

\beginref \parindent=20pt\item{\bref{14}}
\hldest{xyz}{}{bib14}%
Ben Green
and Terence Tao,
``The primes contain arbitrarily long arithmetic progressions,''
{\sl Annals of Mathematics}\/
{\bf 167}
(2008),
481--547.
\endref

\beginref \parindent=20pt\item{\bref{15}}
\hldest{xyz}{}{bib15}%
Curtis Greene,
``A class of lattices with M\"obius function $\pm 1,0$,''
{\sl European Journal of Combinatorics}\/
{\bf 9}
(1988),
225--240.
\endref

\beginref \parindent=20pt\item{\bref{16}}
\hldest{xyz}{}{bib16}%
Becky Hooper
and Adam Van Tuyl,
``A note on the van der Waerden complex,''
{\sl Mathematica Scandinavica}\/
{\bf 124}
(2019),
179--187.
\endref

\beginref \parindent=20pt\item{\bref{17}}
\hldest{xyz}{}{bib17}%
Petteri Kaski,
Jukka Kohonen,
and Thomas Westerb\"ack,
``Fast M\"obius inversion in semimodular lattices and ER-labelable posets,''
{\sl Electronic Journal of Combinatorics}\/
{\bf 23}
(2016),
P3.26.
\endref

\beginref \parindent=20pt\item{\bref{18}}
\hldest{xyz}{}{bib18}%
Tiansi Li,
``EL-shelling on comodernistic lattices,''
{\sl Journal of Combinatorial Theory, Series A}\/
{\bf 177}
(2021).
\endref

\beginref \parindent=20pt\item{\bref{19}}
\hldest{xyz}{}{bib19}%
James Munkres,
{\sl Topology, 2nd edition}
(Upper Saddle River:
Prentice Hall,
2000).
\endref

\beginref \parindent=20pt\item{\bref{20}}
\hldest{xyz}{}{bib20}%
Gian-Carlo Rota,
``On the foundations of combinatorial theory I. Theory of M\"obius Functions,''
{\sl Zeitschrift f\"ur Wahrscheinlichkeitstheorie und verwandte Gebiete}\/
{\bf 2}
(1964),
340--368.
\endref

\beginref \parindent=20pt\item{\bref{21}}
\hldest{xyz}{}{bib21}%
Jay Schweig
and Russ Woodroofe,
``A broad class of shellable lattices,''
{\sl Advances in Mathematics}\/
{\bf 313}
(2017),
537--563.
\endref

\beginref \parindent=20pt\item{\bref{22}}
\hldest{xyz}{}{bib22}%
Endre Szemer\'edi,
``On sets of integers containing no $k$ elements in arithmetic progression,''
{\sl Acta Arithmetica}\/
{\bf 27}
(1975),
199--245.
\endref

\beginref \parindent=20pt\item{\bref{23}}
\hldest{xyz}{}{bib23}%
Bartel Leendert van~der~Waerden,
``Beweis einer Baudetschen Vermutung,''
{\sl Nieuw Archief voor Wiskunde}\/
{\bf 15}
(1927),
212--216.
\endref

\goodbreak\
        \bye

%% file: fontmac.tex


\font\smallheader=cmssbx10 

\font\eightpt=cmr8
\font\ninept=cmr9


\font\ninerm=cmr9     \font\eightrm=cmr8   \font\sixrm=cmr6      
\font\ninei=cmmi9     \font\eighti=cmmi8   \font\sixi=cmmi6      
\font\ninesy=cmsy9    \font\eightsy=cmsy8  \font\sixsy=cmsy6     
\font\ninebf=cmbx9    \font\eightbf=cmbx8  \font\sixbf=cmbx6     
\font\ninett=cmtt9    \font\eighttt=cmtt8                        
\font\nineit=cmti9    \font\eightit=cmti8     
\font\ninesl=cmsl9    \font\eightsl=cmsl8                        

\font\tensc=cmcsc10   \font\ninesc=cmcsc9  \font\eightsc=cmcsc8  

\font\eightssq=cmssq8  \font\eightssqi=cmssqi8  


\font\tenssbx=cmssbx10 

  \font\twelvebf=cmbx12
  
\def\sc{\tensc}  \def\mc{\ninerm}

\input cyracc.def
    \font\tencyr=wncyr10   \font\ninecyr=wncyr9   \font\eightcyr=wncyr8
    \font\tencyri=wncyi10  \font\ninecyri=wncyi9  \font\eightcyri=wncyi8
    \def\cyr{\tencyr\cyracc} \def\cyri{\tencyri\cyracc}

\newskip\ttglue
\def\tenpoint{\def\rm{\fam0\tenrm}%
  \textfont0=\tenrm \scriptfont0=\sevenrm \scriptscriptfont0=\fiverm
  \textfont1=\teni  \scriptfont1=\seveni  \scriptscriptfont1=\fivei
  \textfont2=\tensy \scriptfont2=\sevensy \scriptscriptfont2=\fivesy
  \textfont3=\tenex \scriptfont3=\tenex   \scriptscriptfont3=\tenex
  \textfont\itfam=\tenit  \def\it{\fam\itfam\tenit}%
  \textfont\slfam=\tensl  \def\sl{\fam\slfam\tensl}%
  \textfont\ttfam=\tentt  \def\tt{\fam\ttfam\tentt}%
  \textfont\bffam=\tenbf  \scriptfont\bffam=\sevenbf
   \scriptscriptfont\bffam=\fivebf \def\bf{\fam\bffam\tenbf}%
  \tt \ttglue=.5em plus.25em minus.15em
  \normalbaselineskip=12pt
  \setbox\strutbox=\hbox{\vrule height8.5pt depth3.5pt width0pt}%
  \let\sc=\tensc \let\mc=\ninerm  
  \def\cyr{\tencyr\cyracc}\def\cyri{\tencyri\cyracc}
  \let\big=\tenbig  \normalbaselines\rm}

\def\ninepoint{\def\rm{\fam0\ninerm}%
\textfont0=\ninerm  \scriptfont0=\sixrm  \scriptscriptfont0=\fiverm
\textfont1=\ninei   \scriptfont1=\sixi   \scriptscriptfont1=\fivei
\textfont2=\ninesy  \scriptfont2=\sixsy  \scriptscriptfont2=\fivesy
\textfont3=\tenex   \scriptfont3=\tenex  \scriptscriptfont3=\tenex
\textfont\itfam=\nineit  \def\it{\fam\itfam\nineit}%
\textfont\slfam=\ninesl  \def\sl{\fam\slfam\ninesl}%
\textfont\ttfam=\ninett  \def\tt{\fam\ttfam\ninett}%
\textfont\bffam=\ninebf  \scriptfont\bffam=\sixbf
\scriptscriptfont\bffam=\fivebf\def\bf{\fam\bffam\ninebf}%
\tt\ttglue=.5em plus.25em minus.15em
\normalbaselineskip=11pt
\setbox\strutbox=\hbox{\vrule height8pt depth3pt width0pt}%
\let\sc=\ninesc\let\mc=\eightrm
\def\cyr{\ninecyr\cyracc}\def\cyri{\ninecyri\cyracc}
\let\big=\ninebig\normalbaselines\rm}

\def\eightpoint{\def\rm{\fam0\eightrm}%
  \textfont0=\eightrm \scriptfont0=\sixrm \scriptscriptfont0=\fiverm
  \textfont1=\eighti  \scriptfont1=\sixi  \scriptscriptfont1=\fivei
  \textfont2=\eightsy \scriptfont2=\sixsy \scriptscriptfont2=\fivesy
  \textfont3=\tenex   \scriptfont3=\tenex \scriptscriptfont3=\tenex
  \textfont\itfam=\eightit  \def\it{\fam\itfam\eightit}%
  \textfont\slfam=\eightsl  \def\sl{\fam\slfam\eightsl}%
  \textfont\ttfam=\eighttt  \def\tt{\fam\ttfam\eighttt}%
  \textfont\bffam=\eightbf  \scriptfont\bffam=\sixbf
  \normalbaselineskip=9pt
  \let\sc=\eightsc \let\mc=\sevenrm  
  \def\cyr{\eightcyr\cyracc}\def\cyri{\eightcyri\cyracc}
  \let\big=\eightbig  \normalbaselines\rm}%
\def\nospace{\nulldelimiterspace0pt\mathsurround0pt}%
\def\tenbig#1{{\hbox{$\left#1\vbox to8.5pt{}\right.\nospace$}}}%
\def\ninebig#1{{\hbox{$\textfont0=\tenrm\textfont2=\tensy
  \left#1\vbox to7.25pt{}\right.\nospace$}}}%
\def\eightbig#1{{\hbox{$\textfont0=\ninerm\textfont2=\ninesy
  \left#1\vbox to6.5pt{}\right.\nospace$}}}%

\def\nonextendedbold{
  \font\fiveb=cmb10 at 5pt
  \font\sixb=cmb10 at 6pt
  \font\sevenb=cmb10 at 7pt
  \font\eightb=cmb10 at 8pt
  \font\nineb=cmb10 at 9pt
  \font\tenb=cmb10
  \font\twelveb=cmb10 at 12pt
  \let\fivebf=\fiveb
  \let\sixbf=\sixb
  \let\sevenbf=\sevenb
  \let\eightbf=\eightb
  \let\ninebf=\nineb
  \let\tenbf=\tenb
  \let\twelvebf=\twelveb
}

\def\leftrighttop#1#2{
  \headline{\ifnum\pageno=1\hfil\else{\ninept #1 \hfil #2}\fi}
}

\def\firstnopagenum{
  \footline{\ifnum\pageno=1 \hfil \else \hfil{\rm \number\pageno}\hfil\fi}
}

\def\maketitle#1#2#3#4{
  \centerline {\titlefont #1}
  \medskip
  \centerline {\eightpt #2}
  \medskip
  \centerline {\tensc #3}
  \medskip
  \centerline {\tensc #4}
  \bigskip
}


\outer\def\floattext#1 #2. #3\par{
  $$
  \vbox{
    \hsize #1 true in
    \noindent{\bf #2.}\enskip #3
  }
  $$
}


\def\lsection#1\par{
  \bigskip\vskip\parskip
  \leftline{\sectionfont#1}\nobreak\medskip\noindent
}

\def\csection#1\par{
  \bigskip\vskip\parskip
  \centerline{\sectionfont#1}\nobreak\medskip\noindent
}

\def\rsection#1\par{
  \bigskip\vskip\parskip
  \rightline{\sectionfont#1}\nobreak\medskip\noindent
}
\def\section{\lsection}

\def\boldlabel#1. {\noindent{\bf #1.\enspace}}
\def\subsection#1. {\medskip\noindent{\bf #1.\enspace}}



\font\tenfrak=eufm10
\font\sevenfrak=eufm7
\font\fivefrak=eufm5
\newfam\frakfam
\textfont\frakfam=\tenfrak
\scriptfont\frakfam=\sevenfrak
\scriptscriptfont\frakfam=\fivefrak
\def\frak#1{{\fam\frakfam #1}}

\def\janksc#1#2 {#1{\eightpt#2}}
\def\jankscsp#1#2 {#1{\eightpt#2}\ }
\def\scproclaim#1.#2\par{\noindent\jankscsp #1.\enspace{\it#2\par}}

\def\ref#1{[#1]}

\def\quote{
  \begingroup
    \baselineskip 10pt
    \parfillskip 0pt
    \interlinepenalty 10000
    \leftskip 0pt plus 40pc minus \parindent
    \let\rm=\quoterm\let\sl=\quotesl\everypar{\sl}
    \obeylines
}
\def\author#1(#2){\nobreak\smallskip\rm--- \rm#1\unskip\enspace(#2)\par\endgroup}

\def\titlefont{\twelvebf}
\def\sectionfont{\tenssbx}
\def\quoterm{\eightssq}
\def\quotesl{\eightssqi}


\tenpoint

%% file: mathmac.tex


\def\xskip{\hskip 7pt plus 3pt minus 4pt}

\def\proof{\medbreak\noindent{\it Proof.}\xskip\ignorespaces}

\def\slug{\quad\hbox{\kern1.5pt\vrule width2.5pt height6pt depth1.5pt\kern1.5pt}\medskip}
\def\noskipslug{\quad\hbox{\kern1.5pt\vrule width2.5pt height6pt depth1.5pt\kern1.5pt}}

\newdimen\algindent
\newif\ifitempar \itempartrue 
\def\algindentset#1{\setbox0\hbox{{\bf #1.\kern.25em}}\algindent=\wd0\relax}
\def\algbegin #1 #2{\algindentset{#21}\alg #1 #2} 
\def\aalgbegin #1 #2{\algindentset{#211}\alg #1 #2} 
\def\alg#1(#2). {\medbreak 
  \noindent{\bf#1}({\it#2\/}).\xskip\ignorespaces}
\def\algstep#1.{\ifitempar\smallskip\noindent\else\itempartrue
  \hskip-\parindent\fi
  \hbox to\algindent{\bf\hfil #1.\kern.25em}%
  \hangindent=\algindent\hangafter=1\ignorespaces}


\def\NN{{\bf N}}
\def\ZZ{{\bf Z}}

\def\CC{{\bf C}}

\def\op#1{\mathop{\hbox{#1}}\nolimits}




\newcount\thmcount  
\thmcount=1
\newcount\sectcount  
\sectcount=1
\newcount\figcount  
\figcount=1
\newcount\eqcount  
\eqcount=1

\def\oldno#1{\eqno({\oldstyle#1})}
\def\refeq#1{({\oldstyle#1})}
\def\adveq{\oldno{\the\eqcount}\global\advance\eqcount by 1}  
\def\advthm{\the\thmcount\global\advance \thmcount by 1}

\def\advsect{\section\the\sectcount\global\advance\sectcount by 1. }

\def\caption#1{\centerline{\ninepoint{\bf Fig.~\the\figcount\global\advance\figcount by 1.\enspace}#1}}

\outer\def\parenproclaim #1 (#2).#3\par{\medbreak
  \noindent{\bf #1}\enspace\rm({\it #2\/}).\nobreak\ignorespaces{\sl #3\par}
  \ifdim\lastskip<\medskipamount \removelastskip\penalty55\medskip\fi}


\newdimen\axiomindent
\def\axiomindentset#1{\setbox0\hbox{{\bf #1.\kern.25em}}\axiomindent=\wd0\relax}
\def\axiom#1. [#2.]{\ifitempar\par\noindent\else\itempartrue
  \hskip-\parindent\fi%
  \hbox to\axiomindent{\bf\hfil #1.\kern.25em}%
  \hangindent=\axiomindent\hangafter=1[{\it #2.}]}